\documentclass{amsart}
\usepackage{amssymb}
\usepackage{amsmath}
\usepackage{amsfonts}
\theoremstyle{plain}
\newtheorem{thm}{Theorem}[section]

\newtheorem{lmm}[thm]{Lemma}
\newtheorem{cor}[thm]{Corollary}
\newtheorem{dfn}[thm]{Definition}

\theoremstyle{remark}

\def\II{\mathbb{I}}
\def\pmc#1{\setbox0=\hbox{#1}
    \kern-.1em\copy0\kern-\wd0
    \kern.1em\copy0\kern-\wd0}

\def\om{\omega}
\def\op{\operatorname}
\def\ov{\overline}

\def\sm{\setminus}

\begin{document}
\title
{On the second homotopy group of $SC(Z)$}

\author{Katsuya Eda}
\address{School of Science and Engineering,
Waseda University, Tokyo 169-8555, Japan}
\email{eda@logic.info.waseda.ac.jp}

\author{Umed H. Karimov}
\address{Institute of Mathematics,
Academy of Sciences of Tajikistan,
Ul. Ainy $299^A$, Dushanbe 734063, Tajikistan}
\email{umed-karimov@mail.ru}

\author{Du\v san Repov\v s}
\address{Faculty of Mathematics and Physics, and
Faculty of Education,
University of Ljubljana, P.O.Box 2964,
Ljubljana 1001, Slovenia}
\email{dusan.repovs@guest.arnes.si}
 
\date{\today}



\begin{abstract}
In our earlier paper we introduced a cone-like space $SC(Z)$. In the present note we
establish some new algebraic properties of  $SC(Z)$.
\end{abstract}
\date{\today}
\maketitle 

\section{Introduction}
In our earlier 
paper 
\cite{EKR:topsin} we introduced a new, cone-like construction of a
space $SC(Z)$ using the topologist sine curve and proved that $SC(Z)$ is
simply connected for every path-connected space $Z$.  
In another paper \cite{EKR:nonaspherical} we proved that the singular homology
$H_2(SC(Z);\mathbb Z)$ is non-trivial if $\pi _1(Z)$ is non-trivial. In the
present paper we prove its converse, that is: 
\begin{thm}\label{thm:main}
Let $Z$ be any path-connected space, $z_{0}\in Z$. If $\pi _1(Z,z_{0})$ is trivial, then $\pi
 _2(SC(Z),z_{0})$ is also trivial. 
\end{thm}
Consequently, we get the following:
\begin{cor}\label{cor:main}
For any path-connected space $Z$, $z_{0}\in Z$, the following statements are equivalent:
\begin{itemize}
\item[(1)] $\pi _1(Z,z_{0})$ is trivial; 
\item[(2)] $\pi _2(SC(Z),z_{0})$ is trivial; and
\item[(3)] $H_2(SC(Z);\mathbb Z)$ is trivial.
\end{itemize}
\end{cor}
We also use this opportunity
to correct the proof of Lemma 3.2 in
\cite{EKR:casechamberin} (see Section 3). 

\section{Proof of Theorem~\ref{thm:main}}

In order to describe the homotopies we shall need to introduce some notations. 
The unit interval $[0,1]$ is denoted by $\II$. 
For a map $f:[a_0, b_0]\times [c,d] \to X$, define $f^-:[a_0,b_0]\times
[c,d] \to X$ by: 
\[
 f^-(x,y) = f(a_0+b_0-x,y)
\]
and $f_{[a,b]}:[a,b]\times [c,d] \to X$ by: 
\[
 f_{[a,b]}(x,y) = f(a_0 + (b_0-a_0)(x-a)/(b-a),y). 
\]
We follow the notation in \cite{EKR:topsin} for the space $SC(Z)$ and
the projection $p: SC(Z)\to \II^2$. 
In particular, a point of the subspace $\II ^2$ of $SC(Z)$ is denoted by
$(x; y)$. The following figure denotes
the part $\II^2$ of $SC(Z)$, where the polygonal
line
$A_1B_1A_2B_2\cdots AB$ is the piecewise linear version of the
topologists sine curve in Figure 1.

\begin{dfn}\label{dfn:standard}\cite[Definition 2.1]{EK:aloha}
A continuous map $f:\II^2\to p^{-1}(\II\times \{ 0\})$
 with $f(\partial \II^2)=\{ A\}$ is
 said to be {\it standard,} if 
$f([1/(m+1), 1/m] \times \II )\subset [A,A_m]\cup p^{-1}(\{ A_m\})$,
and $f(\partial ([1/(m+1), 1/m] \times \II)) = \{ A\}$ 
for each $m<\om$. 
\end{dfn}

\begin{center}
\setlength{\unitlength}{0.7mm}
\begin{picture}(110,110)
\put(26,101){\line(-1,-6){10}}
\put(26,0){\line(0,1){101}}
\put(26,0){\line(1,4){25}}
\put(51,0){\line(0,1){101}}
\put(51,0){\line(1,2){50}}
\put(1,1){\line(1,0){100}}
\put(1,101){\line(1,0){100}}
\put(1,1){\line(0,1){100}}
\put(101,1){\line(0,1){100}}

\put(-3,-4){\shortstack{$A$}}
\put(-2,104){\shortstack{$B$}}
\put(102,104){\shortstack{$B_1$}}
\put(50,104){\shortstack{$B_2$}}
\put(25,104){\shortstack{$B_3$}}
\put(102,-4){\shortstack{$A_1$}} 
\put(50,-4){\shortstack{$A_2$}} 
\put(25,-4){\shortstack{$A_3$}} 
\put(100,0){\shortstack{$\bullet$}} 
\put(50,0){\shortstack{$\bullet$}} 
\put(25,0){\shortstack{$\bullet$}} 
\put(25,100){\shortstack{$\bullet$}} 
\put(50,100){\shortstack{$\bullet$}} 
\put(100,100){\shortstack{$\bullet$}}
\put(0,100){\shortstack{$\bullet$}}
\put(0,0){\shortstack{$\bullet$}} 
\end{picture}
\end{center}

\smallskip
\begin{center}
(Figure 1)
\end{center}
\smallskip

{\it Proof of\/} Theorem~\ref{thm:main}.
Recall that
$p: SC(Z)\to \II^2$ and $p(z)= (p_1(z) ; p_2(z))$. 
Let $f: \II^2\to SC(Z)$ with $f(\partial \II^2) = A =
(0,0)$. First we claim that $f$ is homotopic to a map $f_1$ in
$p^{-1}(\II\times \{ 0\})$ relative to $\partial
\II^2$. 

Since $p^{-1}(\II\times \{ 0\})$ is a retraction
of $SC(Z)\sm \bigcup \{ p^{-1}(\{ B_n\}): n\in \mathbb{N}\}$, it
suffices to show that $f$ is 
homotopic to a map in $SC(Z)\sm \bigcup \{ p^{-1}(\{ B_n\}) : n\in
\mathbb{N}\}$ relative to $\partial\II^2$.
The idea to remove
$p^{-1}(\{ B_n\})$ from the image of $f$ is basically the same as
that in the proof of the simple connectivity of $SC(Z)$. 
Since the
boundary of an open connected subset of the square is complicated in
comparison with those of the interval, some more care is needed.

Here we also use the simple connectivity of $Z$. Using this
property we get the simple connectivity of a certain small open subset
containing $p^{-1}(\{ B_n\})$ in $SC(Z)$, i.e. $p^{-1}(U_n)$ according 
to the following
notation. Let $U_n$ be a square neighborhood of $B_n$ in
$\II ^2$ which does not contain any $B_i$ 
for $i\neq n$, and choose a point $u_n\in U_n\cap T$ satisfying $u_n\neq B_n$. 

We have finitely many polygonal connected open
sets $O_i$ in $\II^2$ such that 
$(p\circ f)^{-1}(\{ B_n\})\subseteq
\bigcup _i O_i \subseteq (p\circ f)^{-1}(U_n)$, 
where $O_i$ may fail to 
be simply connected.  
Observe that $p^{-1}(U_n)$ is homotopy equivalent to
$Z$ and so it
is simply connected. 
Hence, working in each $O_i$,
we have $g_1: \II^2\to SC(Z)$ satisfying: 
\begin{itemize}
\item[(1)] $g_1$ is homotopic to $f$ relative to $\II ^2\sm O_i$; and  
\item[(2)] there exist finitely many simply connected polygonal pairwise
	   disjoint subsets $P_{ij}$ of $O_i$ such that
	   $(p\circ g_1)^{-1}(\{ B_n\})) \subseteq \bigcup _j P_{ij}$, and 
\item[(3)] $g(\partial P_{ij}) = \{ u_n\}$. 
\end{itemize}

We remark that the range of the homotopy between $f$ and $g_1$ is
contained in $p^{-1}(U_n$). Since $p^{-1}(\{ u_n\})$ is a strong
deformation retract
of $p^{-1}(U_n)$, we have $g_2: \II^2\to SC(Z)$ such that
$g_2(\bigcup _j P_{ij}) \subseteq p^{-1}(\{ u_n\})\cup \II\times \II$, 
$g_2$ is homotopic $g_1$ relative to $\II^2\sm \bigcup _j P_{ij}$, and
the range of the homotopy is contained in $p^{-1}(U_n)$. 

Observe that there are only finitely many $O_i$ for each $B_n$ and
that $p^{-1}(U_n)$ converge to $B$. Working on each $B_n$ successively, we
obtain maps homotopic to $f$. 
Let $f_0$ be the limit of these maps. 

Since $f(O_i)\subseteq p^{-1}(U_n)$ and $p^{-1}(U_n)$ converge
to $B$, $f_0$ is continuous and $f_0$ is homotopic to $f$ relative to
$\partial\II^2$ and also $f_0(\II^2)$ does not intersect with any
$p^{-1}(\{B_n\})$, as desired. Hence we have 
$f_1: \II^2\to p^{-1}(\II\times \{ 0\})$ which is
homotopic to $f$ relative to $\partial\II^2$. 

The next procedure is similar to the proof of \cite[Lemma 2.2]{EK:aloha},
which is rather long but each 
step
is simple.  
Using the commutativity of $\pi _2$ and the simple connectivity of $Z$
again and also using the fact that $p^{-1}(\II\times \{ 0\})$ is locally
strongly contractible at points $(x;0)$ with $(x;0)\notin \{ A_n: n\in
\mathbb{N}\}\cup \{ A\}$, 
we 
get
a standard map $f_2: \II^2\to p^{-1}(\II\times
\{ 0\})$ which is homotopic to $f_1$ relative to $\partial\II^2$. 

The proof that $f_2$ is null-homotopic is the $2$-dimensional version of
procedures II and III in the proof of the
simple connectivity of $SC(Z)$ \cite[Theorem 1.1]{EKR:topsin}. 
We outline these procedures. We concentrate on description 
of the null homotopy of $f_2\, | \, [1/(k+1),1/k]\times \II$. 

Fix $k\in \mathbb{N}$. For $m\in \mathbb{N}$,
define $h_m: [1/(k+1),1/k]\times \II \to SC(Z)$ by:  
\[
h_m(x,y) = \left\{
\begin{array}{ll}
(ku/(k+m-1); 0) \quad   &\mbox{if } f_2(x,y)=(u; 0)
 \\
(A_{k+m-1},z) \quad   &\mbox{if } f_2(x,y)=(A_k,z).
\end{array}\right.
\]
Next define 
$g_{k,m}: [1/(k+1) + 1/((m+1)k(k+1)),1/(k+1) + 1/(mk(k+1))]\times \II
\to SC(Z)$ by: 
\begin{eqnarray*}
g_{k,2m-1} &=& (h_m)_{[1/(k+1) + 1/(2mk(k+1)),1/(k+1) + 1/((2m-1)k(k+1))]}
 \\
g_{k,2m} &=& (h^- _{m+1})_{[1/(k+1) + 1/((2m+1)k(k+1)),1/(k+1) + 1/(2mk(k+1))]}
\end{eqnarray*}
Let $g_k: [1/(k+1),1/k] \times \II\to SC(Z)$ be the unique
continuous extension of 
\[
\bigcup _{m\in \mathbb{N}} g_ {k,m}, 
\]
i.e. $g_k(1/(k+1),y) = A$. 

Since the images of $g_{k,m}$ converge to $A$ and $g_{2m+1}$ is the
homotopy inverse of $g_{2m}$ in $p^{-1}(\II \times \{ 0\})$ for each
$m$, $g_k$ is continuous and is homotopic relative to the boundary to
the restriction $f_2\, |\, [1/(k+1),1/k]\times \II$,
and the homotopy can
be taken
in $p^{-1}([A,A_k])$. Hence $f_2$ is homotopic relative to the boundary
to $g:\II\times \II\to SC(Z)$ which is the unique continuous extension
of $\bigcup \{ g_k: k\in \mathbb{N}\}$, i.e. $g(0,y) = A$.  

For the next step we do not
care for the boundary for a while. 
We push up the ranges of $g_{k,2m-1}$ along $A_{k+m}B_{k+m}$ for $m\ge
0$ and $g^-_{k,2m}$ along $A_{k+m}B_{k+m-1}$ so that the
$y$-coordinates of $p(u)$ for $u$ in each of the ranges are the same. 
Then the resulting map is defined in 
$p^{-1}(\II\times \{ 1\})$ and we couple $g_{k,1}$ and $g_ {k,2}$, and
generally $g_{k,2m-1}$ and $g_{k,2m}$. 
Since these homotopies of couplings converge to $B$, we see that the
resulting map is null-homotopic. 
We can perform these procedures uniformly in $k$, and we have
a homotopy from $f_2$ to the constant map $B$. 

In order to obtain the desired homotopy to the constant $A$ relative to
the boundary, we can modify the homotopy above to the desired one, because 
we have homotopies in the  pushing up procedure above, so that the
$y$-coordinates of $p(u)$ for $u$ in the ranges are the same, even
uniformly in $k$.
\qed

\section{Correction of the proof of Lemma 3.2 of \cite{EKR:casechamberin}}
 
In our earlier paper \cite[Lemma 3.2]{EKR:casechamberin} 
we used the following auxiliary result:

\begin{lmm}\label{prop:loop}
Let $p_0,p_1,p^*$ be distinct points in a Hausdorff space $X$ and let
$f$ be
 a loop in $X$ with the
 base point $p_1$ such that $f^{-1}(\{ p_0\})=\emptyset$ 
 and $f^{-1}(\{ p^*\})$ is a singleton. 
 
If $f$ is null-homotopic relative to end points, 
then there exists a loop 
 $f'$ in $X$
 with the base point $p_1$ in $X\setminus \{ p_0,p^*\}$ such that $f$
 and $f'$ are homotopic relative to end points in $X\setminus \{ p_0\}$. 
\end{lmm}

We use this opportunity to correct our original proof. The assertion
``$G^{-1}(\{ p^*, p_0\})\cap O$ is compact'' in \cite[p.92 l.5 from the
bottom]{EKR:casechamberin} is wrong. 
We begin by the following
well-known result - see e.g. \cite[p.169]{Kuratowski:topology2}:

\begin{lmm}\label{lmm:wellknown}
Let $X$ be a compact space and $C$ 
a closed component of $X$. Then
 $C$ is the intersection of clopen sets containing $C$. \qed
\end{lmm}

\noindent
{\it Proof of\/} Lemma~\ref{prop:loop}.
Since $f$ is null-homotopic, we have a homotopy 
$F: \II\times \II \to X$ from $f$ to the
 constant mapping to $p_1$, i.e. 
 $$F(s,0) = f(s), F(s,1) = F(0,t) = F(1,t) = p_1 \
 \hbox{for} \ 
 s,t\in \II.$$ 
Let $\{ s_0\}$ be the singleton $f^{-1}(\{ p^*\})$. Let
$M_0$ be the connectedness component of $F^{-1}(\{ p^*\})$ containing
 $(s_0,0)$, and $O$  the connectedness
 component of 
$\II\times \II \setminus M_0$ containing $\II\times \{ 1\} $. 
Define $G: \II\times \II \to {\mathcal P}^*$ by: 
\[ G(s,t) = \left\{
\begin{array}{ll}
F(s,t) & \mbox{if} \ \ (s,t)\in O, \\
p^* & \mbox{otherwise.} 
\end{array}\right. \]
Then $G$ is also a homotopy from $f$ to the constant mapping to $p_1$ and 
$G^{-1}(\{ p_0\})$ is contained in $O$. 
(Observe that $\partial (\II\times \II )\sm \{ (s_0,0)\} \subseteq O$.)

Put $M_1 = (\II\times \II \sm O)\cup \II \times \{ 0\}$. 
Since $G^{-1}(\{ p_0\})\cap O$ is compact and disjoint
 from $M_1$, we have a polygonal neighborhood $U$ of $M_1$ whose closure is
 disjoint from $G^{-1}(\{ p_0\})$ and also $\II\times \{ 1\}$. 
The boundary of $U$ is a piecewise linear arc connecting 
a point in $\{ 0\}\times (0,1)$ and a point in $\{ 1\}\times (0,1)$. 
We want to get a piecewise linear injective path $g:\II\to \II\times \II$
 such that 
 $$\op{Im}(G\circ g) \subseteq X\setminus \{ p_0,p^*\}, \  
 g(0)\in \{ 0\}\times \II,\, \
 \hbox{and} \ g(1)\in \{ 1\}\times \II $$
 and  $\op{Im}(g)$ divides $\II\times\II$ into two components,
 one of which contains $G^{-1}(\{ p_0\})$ and the other 
 contains $M_1$.

If the boundary of $U$ is disjoint from $G^{-1}(\{ p^*\})$,
 then we have such a path $g$ tracing the boundary. Otherwise, let 
$C_0$ be the intersection of the boundary of $U$ and
 $G^{-1}(\{ p^*\})$. Since $M_0$ is contained in a connected component
 of $G^{-1}(\{ p^*\})$ which is disjoint with $C_0$, we have clopen sets
 $C_1$ and $C_2$ in $G^{-1}(\{ p^*\})\cap \ov{U}$ such that $C_1\cap C_2
 = \emptyset$, $C_0\subseteq C_1$ and $M_0\subseteq C_2$ by
 Lemma~\ref{lmm:wellknown}. 
 
 Observe that
 $C_1$ and $C_2$ are closed subset of $\ov{U}$. Then we can see that $\{
 0\}\times (0,1) \cap U$ and $\{ 1\}\times (0,1) \cap U$ belong to the
 same component of $U\sm G^{-1}(\{ p^*\})$. We have a piecewise linear
 arc between a point in $\{ 0\}\times (0,1)$ and a point in $\{
 1\}\times (0,1)$ in $U$ which does not intersect with $G^{-1}(\{
 p^*\})$ and have a path $g$ with 
the required properties. 
We now see that $G\circ g$ is the desired loop $f'$.
\qed 

\section{Acknowledgements}
 This research was supported by the
 Slovenian Research Agency grants P1-0292-0101-04 and
J1-9643-0101-07.
The first author was supported by the
Grant-in-Aid for Scientific research (C) of Japan 
No. 20540097. 
We thank the referee for comments and suggestions.
 
\providecommand{\bysame}{\leavevmode\hbox to3em{\hrulefill}\thinspace}


\end{document}